\documentclass[11pt]{amsart}         
\usepackage{amscd}                   
\usepackage{amssymb} 

\newtheorem{thm}{Theorem}[section]
\newtheorem{lem}[thm]{Lemma}

\newtheorem{cor}[thm]{Corollary}

\newtheorem{prop}[thm]{Proposition}

\theoremstyle{definition}

\renewcommand{\thecase}{}

\newtheorem{rmk}[thm]{Remark} 
\renewcommand{\thestep}{}

\theoremstyle{remark}

\makeatletter
\def\alphenumi{
  \def\theenumi{\alph{enumi}}
  \def\p@enumi{\theenumi}
  \def\labelenumi{(\@alph\c@enumi)}}
\makeatother

\makeatletter
\def\thecase{\@arabic\c@case}
\makeatother

\numberwithin{equation}{section}

\makeatletter
\def\thestep{\@arabic\c@step}
\makeatother

\newenvironment{pf*}[1]{\begin{proof}[#1]}{\end{proof}}

\newcommand\CC{\mathbb{C}}

\newcommand\NN{\mathbb{N}}

\newcommand\RR{\mathbb{R}}

\newcommand{\cov}{\nabla}

\newcommand\half{{\textstyle{\frac{1}{2}}}}

\newcommand\eps{\varepsilon}

\newcommand\su{{\mathfrak{s}\mathfrak{u}}}
\newcommand\fu{{\mathfrak{u}}}

\newcommand\PU{\operatorname{PU}}

\newcommand\SO{\operatorname{SO}}

\newcommand{\8}{\infty}

\newcommand\dist{\operatorname{dist}}

\newcommand\End{\operatorname{End}}

\newcommand\grad{\operatorname{grad}}

\newcommand\Ker{\operatorname{Ker}}

\newcommand\supp{\operatorname{supp}}

\newcommand\spinc{\text{$\text{spin}^c$ }}

\newcommand\sF{{\mathcal{F}}}

\begin{document}
\title[Kato-Yau Inequality and Decay Estimate]
{A Kato-Yau Inequality and Decay Estimate for Eigenspinors} 
\author[Paul M. N. Feehan]{Paul M. N. Feehan}
\address{Department of Mathematics\\
Ohio State University\\
Columbus, OH 43210, U.S.A.}
\email{feehan@math.ohio-state.edu}
\urladdr{http://www.math.ohio-state.edu/$\sim$feehan/} 
\curraddr{Max Planck Institut f\"ur Mathematik\\
Vivatsgasse 7\\
Bonn, D-53111\\
Germany}
\email{feehan@mpim-bonn.mpg.de}
\thanks{The author was supported in part by NSF grant
DMS 9704174 and, through the Institute for Advanced Study (1998-1999), by
NSF grant DMS 9729992} 
\date{This version: August 9, 2000. First version: March 2,
1999. math.DG/9903021.}  

\maketitle


\pagenumbering{arabic}
\section{Introduction}
\label{sec:Intro}
The purpose of this article is to prove a Kato-Yau inequality for harmonic
spinors and a decay estimate for eigenspinors. We also describe some new
applications to gauge theory---specifically, to estimates used when gluing
and ungluing $\PU(2)$ monopoles \cite{FLGeorgia}, \cite{FL3}.

\subsection{Statement of results}
\label{subsec:Statement}
Let $(X,g)$ be an oriented, Riemannian, smooth four-manifold and let
$\Omega=\Omega(x_0,r_0,r_1)$ denote the annulus $B(x_0,r_1)-\bar
B(x_0,r_0)$, where $0<4r_0<r_1<\8$. Denote $r=\dist_g(x,x_0)$.  Consider a
triple $(X,\Omega,g)$ where the Laplacian $\Delta_g$ on $C^\8(X)$
has a $C^\8$ Green kernel $G_g$ and there is a constant $c_0\geq 1$ such
that
\begin{equation}
\label{eq:GreenKernel}
c_0^{-1}r^{-2} \leq G_g(x,x_0) \leq c_0r^{-2},
\quad x\in\Omega(x_0,r_0,r_1),
\end{equation}
and the scalar curvature, $\kappa_g$, of the Levi-Civita connection
$\cov_g$ obeys
\begin{equation}
\label{eq:ScalarCurvature}
|\kappa_g|(x) \leq c_1\eps r^{-2}
\quad\text{and}\quad
|\kappa_g|(x) \leq c_2\left(r_0^2r^{-4} + r_1^{-2}\right),
\quad x\in\Omega(x_0,r_0,r_1).
\end{equation}
Let $(\rho,W)$ be a \spinc structure \cite{SalamonSWBook} on
$X$ with $C^\8$ unitary connection $A_d$ on $\det(W^+)$ and spin
connection on $W=W^+\oplus W^-$, let $E$ be a Hermitian bundle
over $X$, and let $A$ be an $C^\8$ unitary connection on $E$. Assume
that $A_d$ and $A$ are Yang-Mills connections whose curvatures obey, for
some $\eps>0$ to be specified,
\begin{equation}
\label{eq:SpincCurvature}
\|F_{A_d}\|_{L^2(\Omega(x_0,r_0,r_1))} \leq \eps 
\quad\text{and}\quad \|F_A\|_{L^2(\Omega(x_0,r_0,r_1))} \leq \eps.
\end{equation}
Let $D_A$ be the Dirac operator, on sections of $V = W\otimes E$, defined
by $(\rho,g,A_d,A)$. 

\begin{thm}
\label{thm:Main}
Continue the notation and assumptions
of the preceding paragraph. Given positive constants
$c_0$, $c_1$, $c_2$ and an integer $k\geq 0$, there are positive constants
$\eps(c_0,c_1,c_2)\leq 1$ (independent of $k$) and
$c_3(c_0,c_1,c_2,k)$, with the following significance. If
$0<4r_0<r_1<\8$ and $\phi\in L^2(\Omega(x_0,r_0,r_1),V)$ satisfies
\begin{equation}
\label{eq:HarmonicSpinor}
D_A\phi=0 \quad\text{on}\quad\Omega(x_0,r_0,r_1), 
\end{equation}
and conditions \eqref{eq:GreenKernel}, \eqref{eq:ScalarCurvature}, and
\eqref{eq:SpincCurvature} hold, then for all $x\in \Omega(x_0,2r_0,r_1/2)$,
\begin{equation}
\label{eq:Main}
|\cov_A^k\phi|_g(x) 
\leq 
c_3r^{-k}\left({r_0}{r^{-3}} + {r_1^{-2}}\right)
\|\phi\|_{L^2(\Omega(x_0,r_0,r_1),g)}.
\end{equation}
\end{thm}

If $\phi$ is an eigenspinor of $D_A$, with non-zero eigenvalue $\mu$, then
$\phi$ is harmonic with respect to the Dirac operator associated to the
``Friedrich connection'' \cite{Friedrich} defined by the
connection $A$ and eigenvalue $\mu$ (see \S
\ref{subsec:EigenspinorDecay}). Theorem \ref{thm:Main} then leads to the

\begin{cor}
\label{cor:SEMain}
Continue the hypotheses of Theorem \ref{thm:Main}, except that we allow
\begin{equation}
\label{eq:EigenSpinor}
D_A\phi=\mu\phi \quad\text{on}\quad\Omega(x_0,r_0,r_1)
\end{equation}
for some $\mu\in\RR$, instead of equation \eqref{eq:HarmonicSpinor}, and
also assume
\begin{equation}
\label{eq:ExtraRadiusConstraint}
r_1^2\leq\eps.
\end{equation}
Then the remaining assertions of Theorem \ref{thm:Main} hold for the
eigenspinor $\phi$.
\end{cor}

\subsection{Applications}
\label{subsec:RemarkApplic}
Our hypothesis that $A_d$ and $A$ are Yang-Mills connections is stronger
than necessary for many applications (indeed, we shall discuss one such
application in \S \ref{sec:EllipticDirac}): this is just the simplest way
to state our main result. In practice, it is enough that (i) the curvatures
$F_{A_d}$ and $F_A$ obey $C^0$ decay estimates of the shape
\eqref{eq:RadeDecay} and (ii) that the conclusions of Lemma
\ref{lem:BoundaryDecayEstimate} hold with constant independent of
$A_d$ or $A$.  The latter condition is satisfied, for example, if the
curvatures $F_{A_d}$ and $F_A$ obey $C^l$ decay estimates of the shape
\eqref{eq:GenRadeDecay}. A more detailed discussion of some important ways
in which the hypotheses can be relaxed is given in \S
\ref{subsec:MainThmHolds} and \S \ref{sec:EllipticDirac}. 

If $\phi \in L^2(\RR^4,W\otimes E)\cap \Ker D_A$ then it is well-known that
$|\phi|(x) = O(r^{-3})$, for $r\to\8$, so the estimate \eqref{eq:Main} for
the rate of decay of $L^2$ harmonic spinors on $\RR^4$ is sharp
\cite[Equation (3.2.24)]{DK}. (See \cite[\S 3.3.3]{DK} for an explicit
construction of solutions.) However, as we shall explain in \S
\ref{sec:EllipticDirac}, the more interesting applications of Theorem
\ref{thm:Main} and Corollary \ref{cor:SEMain} arise when $\phi$ is an $L^2$
eigenspinor over an annulus, with the given curvature constraints,
rather than all of $\RR^4$ or $S^4$.  Moreover, while standard elliptic
theory would predict an estimate with the general shape of \eqref{eq:Main},
the crucial point is that we know the explicit dependence of the constant
on the radii $r_0$ and $r_1$ and, indirectly via the bounds
\eqref{eq:ScalarCurvature} and \eqref{eq:SpincCurvature}, on the curvatures
of the connections defining the Dirac operator: as we sketch briefly below
and in \S \ref{sec:EllipticDirac}, this is the significant feature of the
estimate \eqref{eq:Main} which allows us to exploit it in gauge-theoretic
applications such as the problem of gluing (and ungluing) $\PU(2)$
monopoles \cite{FLGeorgia}, \cite{FL1}, \cite{FL2b}, \cite{FL3}.

A key difficulty in attempts to directly adapt Taubes' gluing arguments for
anti-self-dual connections \cite{TauStable} to the case of $\PU(2)$
monopoles is the problem of obtaining useful $C^0\cap L^2_2$ estimates for
{\em negative\/} spinors, namely sections $\phi$ of $C^\8(W^-\otimes E)$;
the Bochner formulas \eqref{eq:Bochner} imply that estimates for {\em
positive\/} spinors, or sections of $W^+\otimes E$, are comparatively
straightforward. If $\phi\in \Ker D_A \cap L^2(W^-\otimes E)$, however, the
shape of the Bochner formulas \eqref{eq:Bochner} implies that elementary
methods do not yield $C^0\cap L^2_2$ bounds on $\phi$ which are uniform
with respect to $F_A$ if $A$ ``bubbles'' in the Uhlenbeck sense
\cite{UhlLp}. See \S \ref{sec:EllipticDirac} for a more detailed explanation of
the difficulty.  While standard elliptic theory yields $C^0$ estimates for
$\phi$ on the complement in $X$ of small balls $B(x_i,\delta_i)$, when $F_A$ is
$C^0$ bounded on such a region (but has curvature which bubbles on the
balls $B(x_i,\delta_i)$), these estimates will not necessarily be uniform with
respect to the ball radii: this is a serious problem in gluing contexts, as
one needs estimates which are uniform with respect to the radii $\delta_i$ as
$\delta_i\to 0$. However, we see from Theorem \ref{thm:Main} that the spinor
$\phi$ is $C^0$ bounded on the complement of the ball $B(x_0,2r_0^{1/3})$,
with constant which is {\em independent\/} of $r_0\to 0$ (and also
$r_1\to\8$). Pointwise decay estimates with this uniformity property for
the curvatures of Yang-Mills connections have been derived by
Donaldson \cite{DonApplic}, R\aa de \cite{Rade}, and Groisser-Parker
\cite{GroisserParkerDecay} and it is these decay estimates which motivate
our hypotheses on the curvatures of the connections $\cov_g$, $A_d$, and
$A$; see \S \ref{subsec:YMDecay} for a fuller account.

\subsection{Main ideas in the proof}
Theorem \ref{thm:Main} can be proved in two quite different ways. The
approach taken in the present article uses a pointwise {\em Kato-Yau
inequality\/}, 
\begin{equation}
\label{eq:IntroKatoYau}
|\cov|\phi||^2\leq \left(\frac{m-1}{m}\right)|\cov_A\phi|^2, 
\end{equation}
which we show is obeyed by $D_A$-harmonic spinors over an $m$-dimensional
manifold, coupled with the Bochner formula \eqref{eq:Bochner} for the Dirac
Laplacian $D_A^2$. The Kato-Yau inequality also holds for an eigenspinor of
$D_A$ with eigenvalue $\mu$ if $A$ is replaced by the Friedrich connection
$\tilde A$ defined by $A$ and $\mu$ (see equation
\eqref{eq:FriedrichConnDimm}). The inequality \eqref{eq:IntroKatoYau} then
leads to a differential inequality for $|\phi|$ and hence a decay estimate. 
See \cite{NakajimaYauTrick} for a survey of
such inequalities and applications \cite{BandoKasueNakajima},
\cite{ItohNakajima}, \cite{NakajimaALE}, \cite{SchoenSimonYau},
\cite{Rade}, whose use was pioneered by Yau in his proof of the Calabi
conjecture \cite{Yau}. As pointed out to us by D. Yang \cite{DYangPrivate},
related differential inequalities for vector-valued harmonic functions on
$\RR^{n+1}_+$ were known to Stein \cite[\S VII.3.1]{Stein}. 

After the preprint version of our article \cite{FeehanKatoV3} was
distributed, the preprints \cite{BransonKato} and
\cite{CalderbankGauduchonHerzlich} became available. These articles
describe Kato-Yau inequalities for injectively elliptic first-order linear
differential operators and compute the optimal Kato-Yau constants in a
general setting; the recent article
\cite{CalderbankGauduchonHerzlichSurvey} contains a nice survey of these
results due, independently, to Branson and to Calderbank, Gauduchon, and
Herzlich.

A second (more cumbersome) proof of Theorem \ref{thm:Main}, described in a
preliminary version of this article \cite{FeehanKatoV3}, considers the
Dirac equation on the cylinder $\RR\times S^3$ and hinges on a calculation
of the first eigenvalue, $9/4$, of the square of the Dirac operator on
$S^3$ for the standard metric \cite{Bar}, \cite{Friedrich}, \cite{Hitchin}.

Neither proof of Theorem \ref{thm:Main} requires us to restrict our
attention to dimension four. We restrict to the special case mainly because
this is where our present applications lie and also for expository reasons,
as the proof is easier to follow.

\subsection{Outline of the article}
We begin in \S \ref{sec:DecayHypothesesMain} by recalling the decay
estimates we shall need for Yang-Mills connections (see \S
\ref{subsec:YMDecay}), as well as describing some natural situations where
the hypotheses of Theorem \ref{thm:Main} are known to hold (see \S
\ref{subsec:MainThmHolds}).  In \S \ref{sec:KatoYauInequalities} we show
that harmonic spinors satisfy the pointwise Kato-Yau inequality
\eqref{eq:IntroKatoYau}. As we explain in \S
\ref{sec:DifferentialInequalities}, this leads to a useful differential
inequality for $\Delta|\phi|^{2/3}$ (giving rise to the $O(r^{-3})$ decay
rate), rather than the weaker inequality for $\Delta|\phi|$ which is a
consequence of the standard Kato estimate (and which would only lead to a
$O(r^{-2})$ decay rate and without the essential appearance of the constant
$r_0$ in \eqref{eq:Main}). When integrated twice, this eventually leads to
the decay estimate \eqref{eq:Main}, as we show in \S
\ref{sec:FirstProofMain}, where we present our proof of Theorem
\ref{thm:Main}. In \S \ref{subsec:EigenspinorDecay} we explain how Theorem
\ref{thm:Main} leads, almost immediately, to \ref{cor:SEMain}. Finally, in
\S \ref{sec:EllipticDirac} we describe an application of Theorem
\ref{thm:Main} to the problem of deriving uniform $L^2_2\cap C^0$ estimates
for harmonic spinors when the connection $A$ may bubble.

\subsubsection*{Acknowledgments}
I am grateful to Tom Mrowka for pointing out that calculations of the Dirac
operator spectrum for $S^3$ are well-known and that such calculations can
be found in \cite{Hitchin}. I am also grateful to Hiraku Nakajima for
directing me to the references \cite{Bourguignon} and \cite{NakajimaALE}
and to Deane Yang for describing the early use of Kato-Yau type
differential inequalities in harmonic analysis \cite[\S VII.3.1]{Stein}.
Finally, I would like to thank the Institute for Advanced Study, Princeton,
and the National Science Foundation, for their generous support during the
preparation of this article. 

\section{Decay estimates for Yang-Mills connections and the hypotheses
of Theorem \ref{thm:Main}} 
\label{sec:DecayHypothesesMain}
We review in \S \ref{subsec:YMDecay} the relevant decay estimates of
Donaldson, Groisser-Parker and R\aa de in order to give some context to
Theorem \ref{thm:Main} and its corollaries, to explain the hypotheses on
the curvature of the metric $g$ and curvatures $F_{A_d} $ and $F_A$, and to
introduce some preliminary material we shall need for the remainder of the
article. In \S \ref{subsec:MainThmHolds} we describe some situations of
wider interest where the conditions on $G_g$, $\kappa_g$, $F_{A_d}$, and
$F_A$ are known to hold.

\subsection{Decay estimates for Yang-Mills connections}
\label{subsec:YMDecay}
The simplest case, due to Donaldson and R\aa
de, is when $X$ is $\RR^4$ with its standard metric.

First, recall that the fundamental decay estimates, for a Yang-Mills
connection $A$ over $B(0,r_1)-\{0\}\subset \RR^4$ with $L^2$-small energy,
\begin{equation}
\label{eq:UhlenbeckBallDecay}
|F_A|(x) \leq C([A],r_1)\frac{1}{r^4} \quad 0 < r \leq r_1/2,
\end{equation}
and, by a conformal diffeomorphism, for a Yang-Mills
connection $A$ over $\RR^4-B(0,r_0)$ with $L^2$-small energy,
\begin{equation}
\label{eq:UhlenbeckBallCompDecay}
|F_A|(x) \leq C([A],r_0)\frac{1}{r^4} \quad 2r_0\leq r < \8,
\end{equation}
were first proved by Uhlenbeck \cite[Corollary 4.2]{UhlRem}. However, it is
very useful in applications to gluing and ungluing anti-self-dual
connections to have pointwise estimates of the above shape, but where (i)
the constant $C$ at most depends on the connection $A$ through the $L^2$
norm of its curvature, $F_A$, (ii) the connection $A$ is only known to be
Yang-Mills or anti-self-dual over an annulus $\Omega(r_0,r_1)$ with finite,
positive inner and outer radii, rather than (as above) over a punctured
ball, $B(0,r_1)-\{0\}$, or the complement of a ball, $\RR^4-B(0,R_1)$, and
(iii) the explicit dependence of $C$ on the radii $r_0$, $r_1$ is known.
The earliest such refinement, for anti-self-dual connections, $A$, was due to
Donaldson \cite[Appendix]{DonApplic}, \cite[Proposition 7.3.3]{DK} and
extended by R\aa de \cite[Theorem 1]{Rade} to the case of Yang-Mills
connections. The essential feature common to the estimate \eqref{eq:Main}
above for harmonic spinors and \eqref{eq:RadeDecay},
\eqref{eq:GPAsympFlatDecay}, and \eqref{eq:GPPuncturedBallDecay} below for
Yang-Mills connections, is that the constants on the right (i) at most
depend on the $L^2$ norm of $F_A$ and (ii) the precise dependence on
the radii $r_0$, $r_1$ is known.

\begin{thm}
\label{thm:RadeDecay}
\cite[Proposition 7.3.3]{DK}, \cite[Theorem 1]{Rade} 
There exist positive constants $c$, $\eps$ with the
following significance.  If $0 < 4r_0 < r_1 < \8$, $A$ is a Yang-Mills
connection on $\Omega(r_0,r_1)\subset\RR^4$, with its Euclidean metric, with
$\|F_A\|_{L^2(\Omega(r_0,r_1))} \leq \eps$, and $r=|x|$, then for $x \in
\Omega(2r_0,r_1/2)$,
\begin{equation}
\label{eq:RadeDecay}
|F_A|(x) 
\leq 
c\left(r_0^2r^{-4} + r_1^{-2}\right)
\|F_A\|_{L^2(\Omega(r_0,r_1))}.
\end{equation}
\end{thm}

Though not mentioned explicitly in \cite{DK}, \cite{Rade}, the proof of Theorem
\ref{thm:RadeDecay} extends to give the following more general decay
estimate:

\begin{cor}
\label{cor:RadeDecay}
Continue the hypotheses of Theorem \ref{thm:RadeDecay}. Then for any
integer $k\geq 0$ we have
\begin{equation}
\label{eq:GenRadeDecay}
|\cov_A^kF_A|(x) 
\leq 
c_k r^{-k}\left(r_0^2r^{-4} + r_1^{-2}\right)
\|F_A\|_{L^2(\Omega(r_0,r_1))}.
\end{equation}
\end{cor}

For example, the corollary follows immediately by combining R\aa de's Lemma
2.2 and Theorem $1'$ (the version of his Theorem 1 for a cylinder
$(t_0,t_1)\times S^3$ in place of the annulus $\Omega(r_0,r_1)$).  

There are a couple of standard, useful situations, which we now discuss,
where a metric $g$ approximates a Euclidean metric well enough that
analogues of Theorem \ref{thm:RadeDecay} hold on such Riemannian manifolds,
though some work is involved in order to adapt R\aa de's argument to these
more general cases. The relevant extensions are due to Groisser and Parker.

First, suppose that $(X,g)$ is an {\em asymptotically flat\/} four-manifold
with a single end. Thus, $X$ is a disjoint union $X_0\cup X_\8$, where
$X_0$ is compact and for some $0<R<\8$ there is a diffeomorphism of
$X_\8\cong \RR^4-B(0,R)$, giving coordinates $\{y^i\}$ on $X_\8$ with
respect to which the metric $g$ has the form
$$
g_{ij} = \delta_{ij} + h_{ij},
$$
where, denoting $r=|y|$, 
\begin{equation}
\label{eq:AsympFlatMetric}
r^2|h_{ij}| 
+ r^3|\partial_\alpha h_{ij}| 
+ r^4|\partial_\alpha\partial_\beta h_{ij}| 
\leq
c(g),
\end{equation}
for some positive constant $c(g)$. Let $\Omega(r_0,r_1)$ denote the annulus
$r_0<r<r_1$ in $\RR^4-B(0,R) \cong X_\8$.

\begin{thm}
\label{thm:GPAsympFlatDecay}
\cite[Theorem 1.2]{GroisserParkerDecay}
Let $(X,g)$ be an asymptotically flat four-manifold and let $E$ be a
Hermitian vector bundle over $X$. Then there exist positive constants $c$,
$\eps$, $R$ such that if $4R \leq 4r_0 \leq r_1 < \8$
and $A$ is a Yang-Mills connection on $\Omega(r_0,r_1)$ with
$\|F_A\|_{L^2(\Omega(r_0,r_1))} \leq \eps$, and $r=|y|$, then for $y \in
\Omega(2r_0,r_1/2)$,
\begin{equation}
\label{eq:GPAsympFlatDecay}
|F_A|_g(y) 
\leq 
c\left(r_0^2r^{-4} + r_1^{-2}\right)
\|F_A\|_{L^2(\Omega(r_0,r_1))}.
\end{equation}
\end{thm}

Though not explicitly proved in \cite{GroisserParkerDecay}, the proof of
Corollary \ref{cor:RadeDecay} also yields more general decay estimates in
the situation of Theorem \ref{thm:GPAsympFlatDecay}:
\begin{equation}
\label{eq:GenGPAsympFlatDecay}
|\cov_A^kF_A|_g(y) 
\leq 
c_kr^{-k}\left(r_0^2r^{-4} + r_1^{-2}\right)
\|F_A\|_{L^2(\Omega(r_0,r_1))}.
\end{equation}
Second, suppose $(X,g)$ is four-manifold of {\em bounded geometry\/},
namely positive injectivity radius and Riemannian curvature bounded in
$C^2$.  In \cite{GroisserParkerDecay} Groisser and Parker first prove
Theorem \ref{thm:GPAsympFlatDecay} and then deduce Theorem
\ref{thm:GPPuncturedBallDecay} below as a corollary using a conformal
diffeomorphism from $B(0,\rho)-\{0\}$, with metric $g$ on $B(0,\rho)$ of
bounded geometry, onto $\RR^4-B(0,\rho^{-1})$, endowed with an induced
metric $g_\8$ which is observed to be asymptotically flat.

\begin{thm}
\cite[Theorem 1.1]{GroisserParkerDecay}
\label{thm:GPPuncturedBallDecay}
Let $(X,g)$ be a Riemannian, smooth four-manifold of bounded geometry. Let
$E$ be a Hermitian vector bundle over $X$. Then there exist positive
constants $c$, $\eps$, $\rho$ such that if $0 < 4r_0 < r_1 \leq \rho$ and
$A$ is a Yang-Mills connection on $\Omega(x_0,r_0,r_1)$ with
$\|F_A\|_{L^2(\Omega(x_0,r_0,r_1))} \leq \eps$, and $r = \dist_g(x,x_0)$,
then for $x\in\Omega(x_0,2r_0,r_1/2)$,
\begin{equation}
\label{eq:GPPuncturedBallDecay}
|F_A|_g(x) 
\leq 
c\left(r_0^2r^{-4} + r_1^{-2}\right)
\|F_A\|_{L^2(\Omega(x_0,r_0,r_1))}.
\end{equation}
\end{thm}

Again, in the situation of Theorem \ref{thm:GPPuncturedBallDecay}, one has
stronger decay estimates:
\begin{equation}
\label{eq:GenGPPuncturedBallDecay}
|\cov_A^kF_A|_g(x) 
\leq 
c_kr^{-k}\left(r_0^2r^{-4} + r_1^{-2}\right)
\|F_A\|_{L^2(\Omega(x_0,r_0,r_1))}.
\end{equation}

\subsection{Some applications where the hypotheses of Theorem
\ref{thm:Main} hold} 
\label{subsec:MainThmHolds}
We begin with the constraints on the Riemannian geometry of $(X,g)$. First,
the hypotheses \eqref{eq:GreenKernel} and \eqref{eq:ScalarCurvature}
obviously hold when $(X,g)$ is $\RR^4$ with its Euclidean metric.  Second,
suppose $(X,g)$ is a Riemannian, smooth four-manifold of bounded geometry.
Let $B(x_0,\rho)\subset X$ be a geodesic ball and suppose that $\rho$ is
much less than the injectivity radius of $(X,g)$. Then Lemma
\ref{lem:HarmonicFunction} implies that the Green kernel $G_g$ obeys
\eqref{eq:GreenKernel}. The scalar curvature $\kappa_g$ obeys
\eqref{eq:ScalarCurvature} since, for $x\in B(x_0,\rho)-\{x_0\}$ and $r =
\dist_g(x,x_0)$, 
$$
|\kappa_g|(x) \leq c_1\eps r^{-2}
\quad\text{and}\quad
|\kappa_g|(x) \leq c_2r_1^{-2} \leq c_2(r_0^2r^{-4}+r_1^{-2}),
$$
where we choose $c_1=\|\kappa_g\|_{L^\8(B(x_0,\rho))}$, $\eps=\rho^{2}\ll
1$, and $c_2 = \rho^2\|\kappa_g\|_{L^\8(B(x_0,\rho))}$, with $r_1\leq
\rho$.  Third, suppose $(X,g)$ is an asymptotically flat
four-manifold. Then Lemma \ref{lem:GPHarmonicFunction} now implies that the
Green kernel $G_g$ obeys \eqref{eq:GreenKernel} while our definition
\eqref{eq:AsympFlatMetric} of an asymptotically flat Riemannian
four-manifold ensures that the scalar curvature $\kappa_g$ obeys
\eqref{eq:ScalarCurvature}.

We now turn to the hypotheses on the curvatures $F_{A_d}$ and $F_A$. First,
it is not necessary that $A_d$ and $A$ be Yang-Mills connections, but
rather that their curvatures satisfy the conditions (i) and (ii) described
at the beginning of \S \ref{subsec:RemarkApplic}.  Second, the connections
$A$ of interest in \cite{FL3}, \cite{FL4} are the connection components of
``approximate $\PU(2)$ monopoles'' obtained by splicing anti-self-dual
connections $A_i$ over $S^4$ onto background connections $A_0$ varying in
an Uhlenbeck-compact family. Thus, over small balls
$B(x_i,4\sqrt{\lambda_i})$ in $X$, the connection $A$ obeys the decay
estimate \eqref{eq:GenRadeDecay} because the connections $A_i$ are
anti-self-dual. On the other hand, over the complement of these balls, the
background connection $A_0$ over $X$ obeys a decay estimate of the shape
\eqref{eq:GenRadeDecay} (with $\|F_{A_0}\|_{L^2(\Omega(r_0,r_1))}$ replaced
by $\eps$) because, for $x\in B(x_0,\rho)-\{x_0\}$, we may write
$$
|\cov_{A_0}^k F_{A_0}|_g(x) 
\leq c_k'\eps r^{-k}r_1^{-2}
\leq c_k'\eps r^{-k}(r_0^2r^{-4}+r_1^{-2}),
$$
where $c_k'=\|\cov_{A_0}^k F_{A_0}\|_{L^\8(B(x_0,\rho))}$ and
$\eps=\rho^{k+2}\ll 1$ with $r_1\leq \rho$. In applications such as those
of \cite{FL3}, when $E$ has complex rank two, if the connection $A$ on $E$
corresponds to a unitary connection $A_e$ on $\det(E)$ and an orthogonal
connection $\hat A$ on $\su(E)$, there is no loss in generality if one assumes
that the determinant connections $A_d$ and $A_e$ are Yang-Mills, and so
Theorems \ref{thm:RadeDecay}, \ref{thm:GPAsympFlatDecay}, or
\ref{thm:GPPuncturedBallDecay}, together with their easy corollaries,
guarantee that the constraints on $A_e$ and $A_d$ are satisfied. Thus one
need only ensure that $\hat A$ also obeys the required estimates.

\section{A Kato-Yau inequality for eigenspinors}
\label{sec:KatoYauInequalities}
Recall that any smooth section $\phi$ of a Riemannian vector bundle with
orthogonal connection $A$ satisfies the pointwise {\em Kato inequality\/}
\cite[Equation (6.20)]{FU}:
$$
|\cov|\phi|| \leq |\cov_A \phi|.
$$
If $\phi$ is not arbitrary, but rather satisfies a differential equation,
then the preceding inequality can sometimes be improved to
$$
(1+\delta)|\cov|\phi|| \leq |\cov_A \phi|,
$$
for some positive constant $\delta$. For example, if $\phi = F_A \in
C^\8(\Lambda^2\otimes\fu(E))$ is the curvature of a Yang-Mills connection
$A$ on a Hermitian vector bundle $E$ over a Riemannian four-manifold $X$, then
$F_A$ satisfies the {\em Kato-Yau inequality\/}
$$
|\cov|F_A||^2 \leq \frac{2}{3}|\cov_A F_A|^2.
$$
See \cite[Lemma 3.1]{Rade} for R\aa de's proof when $A$ is Yang-Mills and
\cite[\S 4.2]{ItohNakajima} for Itoh-Nakajima's argument
when $A$ is (anti-)self-dual. A related result
is proved by Bando-Kasue-Nakajima in \cite[Lemma 4.9]{BandoKasueNakajima},
though they attribute the trick to Yau \cite{Bourguignon},
\cite{NakajimaYauTrick}, \cite{NakajimaALE},
\cite{SchoenSimonYau}, \cite{Yau}.

We shall derive a similar Kato-Yau inequality for harmonic spinors.  Given a
Riemannian $m$-manifold $(X,g)$, let $V$ be a complex Hermitian vector
bundle over $X$ with a unitary connection $\cov_A$, and let $\rho:T^*X\to
\End_\CC(V)$ be a linear map such that for all $\alpha\in\Omega^1(X)$,
\begin{itemize}
\item $\rho(\alpha)^\dagger = -\rho(\alpha)$ and 
$\rho(\alpha)^2 = -g(\alpha,\alpha)$, 
\item $[\cov_A,\rho(\alpha)] = \rho(\cov\alpha)$, where $\cov$ is the
Levi-Civita connection on $T^*X$.
\end{itemize}
The composition $D_A = \rho\circ\cov_A$ of Clifford multiplication
$\rho:T^*X\to\End_\CC(V)$ and the covariant derivative $\cov_A:C^\8(V)\to
C^\8(T^*X\otimes V)$ is a generalized Dirac operator and $(\rho,V)$ is a
(complex) Dirac bundle in the sense of \cite[Definition
II.5.2]{LM}. Now suppose that $\phi\in C^\8(V)$ is $D_A$-{\em harmonic\/}:
\begin{equation}
\label{eq:Dirac}
D_A\phi = \sum_{i=1}^m\rho(e^i)\cov_{A,e_i}\phi = 0.
\end{equation}
Recall that for $f\in C^\8(X)$, we have $\cov f = \langle\,\cdot\,,\grad
f\rangle$, so $|\cov f| = |\grad f|$ and $\cov_e f = \langle e,\grad
f\rangle$ for $e\in C^\8(TX)$. If we choose $e = (|\grad f|)^{-1}\grad f$
at points where $(\cov f)(x) \neq 0$ then we obtain the familiar identity
$$
|\cov_e f| = |\cov f|.
$$
We may suppose without loss of generality that $x\in X$ is a point for
which $|\phi|(x)\neq 0$ and $\cov|\phi|(x) \neq 0$.  
At any such point $x$ we can find an
orthonormal frame $\{e_i\}$ for $TX$ and dual coframe $\{e^i\}$ for
$T^*X$ such that
\begin{equation}
\label{eq:Liftoff}
|\cov|\phi|| = |\cov_{e_1}|\phi|| \leq |\cov_{A,e_1}\phi|.
\end{equation}
Indeed, we can take $e_1 = (|\grad|\phi||)^{-1}\grad|\phi|$ and complete
this to give a local orthonormal frame for $TX$ near $x$. Note that
$|\rho(e')\phi| = |\phi|$ for any $e'\in C^\8(T^*X)$ with $|e'|=1$, so
\begin{align*}
|\cov_{A,e_1}\phi|^2
&=
|\rho(e^1)\cov_{A,e_1}\phi|^2
\leq
\left|\sum_{i=2}^m\rho(e^i)\cov_{A,e_i}\phi\right|^2
\\
&\leq
\left(\sum_{i=2}^m|\cov_{A,e_i}\phi|\right)^2
\leq
(m-1)\sum_{i=2}^m|\cov_{A,e_i}\phi|^2,
\end{align*}
using $(\sum_{i=1}^n a_i)^2 \leq n\sum_{i=1}^n a_i^2$ to obtain the final
inequality. Thus,
\begin{equation}
\label{eq:GettingThere}
m|\cov_{A,e_1}\phi|^2
\leq
(m-1)\sum_{i=1}^m|\cov_{A,e_i}\phi|^2
=
(m-1)|\cov_A\phi|^2.
\end{equation}
Combining inequalities \eqref{eq:Liftoff} and \eqref{eq:GettingThere}
yields the following Kato-Yau type inequality on the open set where
$\phi(x)\neq 0$ and $\cov|\phi|(x)\neq 0$. The inequality trivially holds
where $\cov|\phi|(x)= 0$.

\begin{lem}
\label{lem:GKY}
Let $(\rho,V)$ be a \spinc structure over a Riemannian $m$-manifold $(X,g)$
with spin connection $\cov_A$ on $V$. Then, for any smooth section $\phi
\in C^\8(X,V)$ for which $D_A\phi=0$, the following pointwise inequality
holds on the open subset $\{\phi\neq 0\}\subset X$:
$$
|\cov|\phi||^2
\leq
\left(\frac{m-1}{m}\right)|\cov_A\phi|^2
\quad\text{on $X$}.
$$
\end{lem}

One can extend the preceding inequality to the case of eigenspinors with
non-zero eigenvalue (see, for example, the proof of Theorem 4 in
\cite{CalderbankGauduchonHerzlich}). 
Suppose $D_A\phi=\mu\phi$ for some $\mu\in\RR$. The
(metric) Friedrich connection \cite{Friedrich} associated to $\cov_A$ and
$\mu$ is defined by
\begin{equation}
\label{eq:FriedrichConnDimm}
\cov_{\tilde A,\eta} = \cov_{A,\eta} + \frac{\mu}{m}\rho(g(\cdot,\eta)),
\quad
\eta\in C^\8(TX).
\end{equation}
Then $\phi$ is harmonic with respect to the Dirac operator
$D_{\tilde A}=\rho\circ\cov_{\tilde A}=D_A-\mu$ and so one obtains a Kato-Yau
inequality for eigenspinors:

\begin{lem}
\label{lem:GKYEigenspinor}
Continue the hypotheses of Lemma \ref{lem:GKY}, but suppose more generally
that $(D_A-\mu)\phi=0$ for some $\mu\in\RR$ and let $\tilde A$ be the
Friedrich connection associated with $(A,\mu)$. Then, the following pointwise
inequality holds on the open subset $\{\phi\neq 0\}\subset X$:
$$
|\cov|\phi||^2
\leq
\left(\frac{m-1}{m}\right)|\cov_{\tilde A}\phi|^2
\quad\text{on $X$}.
$$
\end{lem}

\section{Differential inequalities for eigenspinors}
\label{sec:DifferentialInequalities}
Using our Kato-Yau inequality we derive the differential inequalities
satisfied by suitable powers of pointwise norms of eigenspinors.

\subsection{Differential inequalities implied by the standard Kato estimate}
Let $(\rho,W)$ be a \spinc structure on the Riemannian $m$-manifold $(X,g)$
with spin connection on a complex Hermitian bundle $W$ of rank $2^n$, where
$m=2n$ or $2n+1$, let $E$ be a complex Hermitian bundle over $X$
equipped with a unitary connection $A$, and let $V= W\otimes E$.

If $m$ is even (respectively, odd), let $A_d$ denote the fixed unitary
connection on $\det(W^+)$ (respectively, $\det(W)$), where $W = W^+\oplus
W^-$, and give $W$ the spin connection induced by $A_d$ and the Levi-Civita
connection on $T^*X$.  To appreciate the significance of Lemma
\ref{lem:GKY}, note that for any $\phi\in C^\8(V)$ the standard identity
\begin{equation}
\label{eq:FUIdentity}
\Delta|\phi|^2 
=
2\langle\cov_A^*\cov_A\phi,\phi\rangle - 2|\cov_A\phi|^2
\end{equation}
and the usual pointwise Kato inequality \cite[Equation (6.20)]{FU},
\begin{equation}
\label{eq:StandardKato}
|\cov|\phi|| \leq |\cov_A \phi|,
\end{equation}
yields the differential inequality
\cite[Equation (6.21)]{FU} on the subset $\{\phi\neq 0\}\subset X$:
\begin{equation}
\label{eq:FUIdentity2}
\Delta|\phi| \leq |\phi|^{-1}\langle\cov_A^*\cov_A\phi,\phi\rangle.
\end{equation}
Now the Bochner formula for $\cov_A^*\cov_A$ on $C^\8(V)$ takes the
general form \cite[Lemma 4.1]{FL1}
\begin{equation}
\label{eq:Bochner}
D_A^2\phi 
=
\cov_A^*\cov_A\phi + \frac{\kappa_g}{4}\phi + \rho(F_A)\phi  
+ \frac{1}{2}\rho(F_{A_d})\phi,
\end{equation}
where $\kappa_g$ is the scalar curvature of the Levi-Civita connection on
$T^*X$.  (When restricted to $C^\8(V^\pm)$ and $X$ has dimension four,
the curvature terms $F_A$ and $F_{A_d}$ above can be replaced by $F_A^\pm$
and $F_{A_d}^\pm$: this leads to useful, global pointwise estimates for
sections of $V^+$ but not for $V^-$, as $F_A^-$ will not be uniformly $L^p$
bounded with respect to $A$ when $p>2$ in applications of interest
\cite{FL3}.)  When $D_A\phi=0$, the inequality \eqref{eq:FUIdentity2} and
the identity \eqref{eq:Bochner} imply
\begin{equation} 
\label{eq:RawHarmonicSectionDiffIneq}
\Delta|\phi| \leq c(|\kappa_g|+|F_A|+|F_{A_d}|)|\phi|.
\end{equation}
We would therefore only expect the standard Kato inequality to at most
imply an $r^{-2}$ decay estimate for $|\phi|(x)$ on $\RR^4$ with its
standard metric, as we can see from our proof of Theorem \ref{thm:Main},
while our Kato-Yau inequality (Lemma \ref{lem:GKY}) implies a $r^{-3}$
decay estimate for $|\phi|(x)$ as we explain in the next subsection.

\subsection{Differential inequalities implied by the Kato-Yau
estimate} 
\label{subsec:KatoYauAndDiffIneq}
The $r^{-4}$ decay estimates of Groisser-Parker, R\aa de, and Uhlenbeck
hinge on differential inequalities for $\Delta|F_A|^{1/2}$ rather than
$\Delta|F_A|$.  To obtain differential inequality in our case which yields
a decay rate strictly faster than $r^{-2}$, suppose $0<\alpha<1$ and
observe that
\begin{align*}
\Delta|\phi|^{\alpha}
&= 
\alpha|\phi|^{\alpha-2}\left((1-\alpha)(\cov|\phi|)^2 
+ |\phi|\Delta|\phi|\right)
\\
&=
\alpha|\phi|^{\alpha-2}\left((2-\alpha)(\cov|\phi|)^2 
+ \frac{1}{2}\Delta|\phi|^2\right).
\end{align*}
(Our convention, when $X=\RR^4$ with its standard metric,
gives $\Delta = -\sum_{i=1}^4\partial^2/\partial x_i^2$, which has sign
opposite to that used in \cite{GT}.) Combining
this with identity \eqref{eq:FUIdentity} implies that
$$
\Delta|\phi|^{\alpha}
=
\alpha|\phi|^{\alpha-2}
\left((2-\alpha)(\cov|\phi|)^2 + \langle \cov_A^*\cov_A\phi,\phi\rangle
- |\cov_A\phi|^2\right).
$$
Thus, we see that the standard Kato inequality does not lead to a useful
differential inequality for $\Delta|\phi|^{\alpha}$ when $\alpha<1$, but
when $D_A\phi=0$, as we now assume, the refinement in Lemma \ref{lem:GKY}
gives $|\cov|\phi||^2 \leq \frac{3}{4}|\cov_A\phi|^2$ when $m=4$ and 
\begin{align*}
\Delta|\phi|^{\alpha}
&\leq
\alpha|\phi|^{\alpha-2}
\left((2-\alpha)\frac{3}{4}|\cov_A\phi|^2
+ \langle \cov_A^*\cov_A\phi,\phi\rangle
- |\cov_A\phi|^2\right)
\\
&=
\alpha|\phi|^{\alpha-2}
\left(((2/3)-\alpha)\frac{3}{4}|\cov_A\phi|^2
+ \langle \cov_A^*\cov_A\phi,\phi\rangle\right).
\end{align*}
Therefore, choosing $\alpha = 2/3$, we obtain
$$
\Delta|\phi|^\alpha 
\leq 
\alpha|\phi|^{\alpha-2}\langle \cov_A^*\cov_A\phi,\phi\rangle.
$$
We now combine the preceding inequality with the Bochner
identity \eqref{eq:Bochner} to give
\begin{equation}
\label{eq:PreRawHarmonicSectionTwoThirdDiffIneq}
\Delta|\phi|^{2/3}
\leq
c(|\kappa_g| + |F_A| + |F_{A_d}|)|\phi|^{2/3} \quad\text{on }X.
\end{equation}
This inequality yields the estimate
\eqref{eq:RawHarmonicSectionTwoThirdDiffIneq} when $\mu=0$. If
$D_A\phi=\mu\phi$, for some $\mu\in\RR$, then estimate
\eqref{eq:PreRawHarmonicSectionTwoThirdDiffIneq} holds with the connection
$A$ replaced by the Friedrich connection \eqref{eq:FriedrichConnDimm}: 
\begin{equation}
\label{eq:FriedrichConnDim4}
\tilde A = A + \frac{\mu}{4}\sum_{i=1}^4\rho(e^i)e^i.
\end{equation}
Thus modified, inequality \eqref{eq:PreRawHarmonicSectionTwoThirdDiffIneq}
yields the estimate \eqref{eq:RawHarmonicSectionTwoThirdDiffIneq} when
$\mu\neq 0$ by writing $F_{\tilde A}$ in terms of $F_A$ and $\mu$. We have
proved:

\begin{lem}
\label{lem:KatoYauAndDiffIneq}
Let $(\rho,W)$ be a \spinc structure over a Riemannian four-manifold
$(X,g)$ with unitary connection $A_d$ on $\det(W^+)$ and spin connection on
the Hermitian rank-$4$ bundle $W=W^+\oplus W^-$ induced by $A_d$ and the
Levi-Civita connection on $T^*X$. Let $E$ be a Hermitian bundle over $X$,
endowed with a unitary connection $A$, let $V^\pm = W^\pm\otimes E$, and
let $V = V^+\oplus V^-$. Let $\kappa_g$ be the scalar curvature of the
Levi-Civita connection. If $\phi\in C^\8(X,V)$ satisfies $D_A\phi=\mu\phi$ on
$X$, for some $\mu\in\RR$, then
\begin{equation} 
\label{eq:RawHarmonicSectionTwoThirdDiffIneq}
\Delta|\phi|^{2/3}
\leq
c(\mu^2 + |\kappa_g| + |F_A| + |F_{A_d}|)|\phi|^{2/3} \quad\text{on }X.
\end{equation} 
\end{lem}

In applications such as those discussed in \cite{FL1}, \cite{FL2b},
\cite{FLGeorgia}, the choice of unitary connection $A_d$ on $\det(W^+)$ is
arbitrary. Typically, therefore, we would assume without loss 
of generality that $A_d$ is
Yang-Mills and so $F_{A_d}$ satisfies the same decay estimates as $F_A$.
{\em For the sake of exposition in the remainder of the article\/}, as it
does not affect the proof, we shall assume $c_1(W^+)=0$ and that $A_d$ is
flat, so $F_{A_d}=0$ in inequalities such as
\eqref{eq:RawHarmonicSectionTwoThirdDiffIneq} above.

As we shall shortly see, the Kato-Yau estimate (and the differential
inequality \eqref{eq:RawHarmonicSectionDiffIneq} which it implies) is the
key ingredient which leads to the stronger $r^{-3}$ estimate of Theorem
\ref{thm:Main}.

As pointed out to us by D. Yang \cite{DYangPrivate}, arguments of the kind
described above leading to useful differential inequalities for
$\Delta|u|^\alpha$ have been known for some time in harmonic analysis, where
$u\in C^\8(\Omega,\RR^{n+1}\otimes \RR^N)$ is a solution to the
generalized Cauchy-Riemann equations on an open subset $\Omega$ of the
half-space $\RR^{n+1}_+$ \cite[\S VII.3.1]{Stein}.

\section{Decay estimates for $L^2$ eigenspinors on Riemannian manifolds}
\label{sec:FirstProofMain}
In this section we prove Theorem \ref{thm:Main} and Corollary \ref{cor:SEMain},
using our differential inequality for eigenspinors (Lemma
\ref{lem:KatoYauAndDiffIneq}).  While our argument  
superficially follows the broad pattern of the proofs of \cite[Theorem
1]{Rade} and \cite[Theorems 1.1 \& 1.2]{GroisserParkerDecay}, the main
new difficulty lies in finding suitable comparison functions for our
applications of the maximum principle, which differ in subtle ways from
those of \cite{GroisserParkerDecay} and \cite{Rade}, as well as a
realization of the correct decay conditions one needs to impose on the
curvatures $F_A$ and $\kappa_g$.

\subsection{Decay estimates for harmonic spinors and proof of Theorem \ref{thm:Main}}
\label{subsec:HarmonicSpinorDecay}
We first deal with the case when the spinor is harmonic.  We begin by
recalling an elementary comparison lemma:

\begin{lem}
\label{lem:BasicComparison}
\cite[Lemma 3.2]{Rade}, 
\cite[Lemma 3.2]{GroisserParkerDecay}
Let $h$ be a positive harmonic function on a domain $\Omega$ in a
Riemannian manifold. Set $\xi = |d(\log\sqrt{h})|^2$. Then for any $a \geq
-1$ the operator $L = \Delta+a\xi$ satisfies the comparison principle on
$\Omega$. Specifically, if $ u, w \in L^2_1(\Omega)$ with $u\leq w$ weakly
on $\partial\Omega$ and $Lu \leq Lw$ weakly on $\Omega$, then $u \leq w$
a.e. on $\Omega$.
\end{lem}

As R\aa de points out in \cite{Rade}, an interesting feature of Lemma
\ref{lem:BasicComparison} is that we are allowed to choose $a<0$ and so
have $a\xi < 0$ on $\Omega$. However, while the comparison
theorem holds for $\Delta+c$, where $c\geq 0$ is any function, it does not
hold in general for $c<0$ on $\Omega$ (for example, a negative constant)
(see \cite[Theorem 3.3]{GT} and the remarks preceding its statement, noting
that their sign conventions for the Laplacian are opposite to ours). 

To apply Lemma \ref{lem:BasicComparison} and ensure that the hypothesis
\eqref{eq:GreenKernel} on the Green kernel of $(X,g)$ is obeyed in the
three geometric situations discussed in \S \ref{subsec:MainThmHolds}, we
need a positive harmonic function $h$ where both $h$ and $\xi =
|d(\log\sqrt{h})|^2$ are essentially $1/r^2$ near the point $x_0$. The case
$X=\RR^4$ with its Euclidean metric is trivial, so we consider the
remaining two situations of interest to us. The first and more subtle
construction, where $(X,g)$ is asymptotically flat, is due to Groisser and
Parker:

\begin{lem}
\label{lem:GPHarmonicFunction}
\cite[Proposition 3.3]{GroisserParkerDecay} Let $(X,g)$ be an
asymptotically flat Riemannian four-manifold. Then there is a positive
constant $\rho_0$ such that for each $\rho<\rho_0$ there is a positive
harmonic function $h$ on $X-B(x_0,\rho^{-1})$ with $h=\rho^2$ on $\partial
B(x_0,\rho^{-1})$ and, if $r = \dist_g(x,x_0)\geq \rho^{-1}$,
\begin{equation}
\label{eq:GPHarmonicFunction}
\frac{3}{4r^2} \leq h \leq \frac{4}{3r^2}
\quad\text{and}\quad
h \leq 8\xi.
\end{equation}
\end{lem}

Second, we consider a small ball in a Riemannian four-manifold $(X,g)$:

\begin{lem}
\label{lem:HarmonicFunction}
Let $(X,g)$ be a Riemannian four-manifold of bounded geometry. Then there
is a positive constant $\rho_0$ such that for each $\rho<\rho_0$ there is a
positive harmonic function $h$ on $B(x_0,\rho)$ with $h=1/\rho^2$ on
$\partial B(x_0,\rho)-\{x_0\}$ and, if $r = \dist_g(x,x_0)\leq\rho$,
\begin{equation}
\label{eq:HarmonicFunction}
\frac{3}{4r^2} \leq h \leq \frac{4}{3r^2}
\quad\text{and}\quad
h \leq 8\xi.
\end{equation}
\end{lem}

\begin{proof}
First, on $\RR^4$ with its standard metric, the Green's function for the
Laplacian on functions is $G(x,y) = -(8\pi^2)^{-1}|x-y|^{-2}$, so we can
choose 
$$
h(x) = \frac{1}{|x|^2} = -8\pi^2G(x,0), \quad x\in\RR^4-\{0\}.
$$
Then $h = \xi$ and so satisfies the constraints in this case. In general,
note that $g_{ij} = \delta_{ij} + O(r^2)$,  $\partial_k
g_{ij} = O(r)$,  and $\partial_k\partial_l g_{ij} =
O(1)$. The existence of $h$ on $B(x_0,\rho)-\{x_0\}$ now follows from its
existence when $g$ is flat and the construction of the Green kernel in
\cite[\S 4.2]{Aubin} for arbitrary metrics.
\end{proof}

\begin{rmk}
The values of the constants $c=3/4$, $c^{-1}=4/3$ in inequalities
\eqref{eq:HarmonicFunction} are not critical. Though any $c \in (0,1)$ would
suffice, the construction of suitable sub- and sup-solutions is a
little easier to follow by the provision of specific constants.
\end{rmk}

For the remainder of this section we shall assume, for the sake of
exposition, that we are in the situation of Lemma
\ref{lem:HarmonicFunction}; there is no essential difference in the proof
for the asymptotically flat case.  Now fix $4r_0<r_1\leq\rho_0$ where
$\rho_0$ is small enough that Lemma \ref{lem:HarmonicFunction} applies. Set
$r(x) = \dist_g(x,x_0)$ on $X$.

In our applications of Lemma \ref{lem:BasicComparison} we shall need an
estimate for $|\phi|$ on the boundaries of the annulus
$\Omega(x_0,2r_0,r_1/2)$, corresponding to the boundary estimates for
$|F_A|$ in equation (3.9) in \cite{GroisserParkerDecay},
where such estimates follow from \cite[Equation
(3.4)]{GroisserParkerDecay}.

\begin{lem}
\label{lem:BoundaryDecayEstimate}
Continue the hypotheses of Theorem \ref{thm:Main}. Then for any integer
$k\geq 0$, there is a constant $c_k$ (independent of the connections $A_d$
on $\det(W^+)$ and $A$ on $E$) such that 
\begin{equation}
\label{eq:BoundaryDecayEstimate}
|\cov_A^k\phi|_g(x) \leq c_k r^{-k-2}\|\phi\|_{L^2(\Omega_g(x_0,r/2,2r),g)}
\quad\text{on }\Omega_g(x_0,2r_0,r_1/2).
\end{equation}
\end{lem}

\begin{proof}
Let $r\in [2r_0,r_1/2]$ be a constant and let $\tilde g = r^{-2}g$, so that
$\Omega_g(r/2,2r) = \Omega_{\tilde g}(1/2, 2)$. Then $D_A^{\tilde g}\phi =
r D_A^g\phi = 0$ and so elliptic estimates on the rescaled annulus
$\Omega_{\tilde g}(1/2, 2)$ and the Sobolev embedding $L^2_3\subset C^0$ give 
\begin{align*}
r^k|\cov_A^k\phi|_{g}(x) 
&= |\cov_A^k\phi|_{\tilde g}(x) 
\\
&\leq
c\|\phi\|_{L^2_{k+3,A}(\Omega_{\tilde g}(3/4, 4/3))}
\quad\text{(on $\Omega_{\tilde g}(3/4, 4/3)$)}
\\
&\leq
c_k\|\phi\|_{L^2(\Omega_{\tilde g}(1/2, 2))}
\\
&=
c_kr^{-2}\|\phi\|_{L^2(\Omega_g(r/2, 2r))},
\end{align*}
and the result follows. The hypotheses of Theorem \ref{thm:Main} on the
connection $A$ ensure that the Sobolev and elliptic estimates above are
uniform with respect to $A$.
\end{proof}

At this point we shall employ the scaling argument in the proof of Lemma
\ref{lem:BoundaryDecayEstimate} to make a further simplification to the
proof of Theorem \ref{thm:Main}.

\begin{lem}
\label{lem:Scaling}
If the estimate \eqref{eq:Main} is valid for some fixed constant $\tilde
r_0>0$ and any $\tilde r_1>4\tilde r_0$, then it holds for any pair of
rescaled constants $(r_0,r_1) = \lambda(\tilde r_0,\tilde r_1)$,
$\lambda>0$, for which the hypotheses of Theorem \ref{thm:Main} are
satisfied.
\end{lem}

\begin{proof}
First observe that the conditions \eqref{eq:GreenKernel},
\eqref{eq:ScalarCurvature}, \eqref{eq:SpincCurvature}, and
\eqref{eq:HarmonicSpinor} are
scale-invariant, as can be seen by replacing the metric $g$ by a rescaled
metric $\tilde g = \lambda^{-2}g$ for any $\lambda > 0$, and using the
rescaling argument of the next paragraph. Hence, the constants $c_0$,
$c_1$, and $c_2$ are also invariant under rescaling.

Let $0<4r_0<r_1<\8$ and define $\lambda>0$ by requiring that $r_0=\lambda
\tilde r_0$. By hypothesis, Theorem \ref{thm:Main} holds for the metric
$\tilde g = \lambda^{-2}g$ and radii $\tilde r_0 = r_0/\lambda$, $\tilde
r_1 = r_1/\lambda >4\tilde r_0$, and $\tilde r = r/\lambda \in [2\tilde
r_0,\tilde r_1/2]$, where $r = \dist_g(x,x_0)$ and $\tilde r =
\dist_{\tilde g}(x,x_0)$. Thus, for $2\tilde r_0 \leq \tilde r \leq \tilde
r_1/2$,
$$
|\cov_A^k\phi|_{\tilde g}(x) 
\leq 
c\tilde r^{-k}
\left(\frac{\tilde r_0}{\tilde r^3} + \frac{1}{\tilde r_1^2}\right)
\|\phi\|_{L^2(\Omega,\tilde g)}
=
cr^{-k}\lambda^k
\left(\frac{\lambda^3 r_0}{\lambda r^3} + \frac{\lambda^2}{r_1^2}\right)
\|\phi\|_{L^2(\Omega,\tilde g)}.
$$
Since $\|\phi\|_{L^2(\Omega,\tilde g)} 
= \lambda^{-2}\|\phi\|_{L^2(\Omega,g)}$ and $|\cov_A^k\phi|_{\tilde g} =
\lambda^k|\cov_A^k\phi|_g$, we see that 
$$
|\cov_A^k\phi|_g(x) 
\leq 
c\left(\frac{r_0}{r^3} + \frac{1}{r_1^2}\right)
\|\phi\|_{L^2(\Omega,g)},
\quad
2r_0 \leq r \leq r_1/2.
$$
as desired.
\end{proof}

Observe that because $r_0<r<r_1$, our hypotheses \eqref{eq:ScalarCurvature}
on the scalar curvature $\kappa_g$ and \eqref{eq:SpincCurvature} on the
curvature $F_A$ (together with the decay estimates \eqref{eq:RadeDecay},
\eqref{eq:GPAsympFlatDecay}, or \eqref{eq:GPPuncturedBallDecay}) imply that
for $r_0 \leq r \leq r_1$,
\begin{align}
\label{eq:RawQuadraticFDecay}
|F_A|_g(x) &\leq c\eps(r_0^2r^{-4}+r_1^{-2}) \leq c\eps r^{-2}, 
\\
\label{eq:RawQuadraticSDecay}
|\kappa_g|(x) &\leq c\eps r^{-2}.
\end{align}
Substituting the $r^{-2}$ decay estimates \eqref{eq:RawQuadraticFDecay} and
\eqref{eq:RawQuadraticSDecay} for $|F_A|$ and 
$\kappa_g$ on $\Omega(x_0,r_0,r_1)$
into the differential inequality \eqref{eq:RawHarmonicSectionDiffIneq}, we
obtain
\begin{equation} 
\label{eq:RiemQuadraticLaplacePhiTwoThird}
\Delta|\phi|^{2/3}
\leq
c\eps r^{-2}|\phi|^{2/3} \quad\text{on }\Omega(x_0,r_0,r_1).
\end{equation} 
Set $E = \|\phi\|_{L^2(\Omega(x_0,r_0,r_1))}$ and fix an arbitrary $\delta
\in (0,1)$.  At this point it simplifies matters considerably if we take
advantage of Lemma \ref{lem:Scaling} and for the remainder of this section,
without loss of generality, assume $r_0\geq 1$. By choosing an $\eps$ small
enough that 
$$
12c\eps\leq \delta(2-\delta),
$$
we see that inequalities
\eqref{eq:HarmonicFunction},
\eqref{eq:BoundaryDecayEstimate}, and
\eqref{eq:RiemQuadraticLaplacePhiTwoThird} imply
\begin{equation}
\label{eq:DiffIneqLaplacePhiTwoThird}
\begin{cases}
|\phi|^{2/3} \leq cE^{2/3}r_0^{-4/3} &\text{on }\partial B(x_0,2r_0),
\\
|\phi|^{2/3} \leq cE^{2/3}r_1^{-4/3} &\text{on }\partial B(x_0,r_1/2),
\\
(\Delta-\delta(2-\delta)\xi)|\phi|^{2/3}
\leq 0 &\text{on } \Omega(x_0,2r_0,r_1/2).
\end{cases}
\end{equation}
(The constraint $\delta < 1$ ensures that $-\delta(2-\delta) > -1$ and
so Lemma \ref{lem:BasicComparison} applies to the operator $L =
\Delta-\delta(2-\delta)\xi$.)

\begin{lem}
\label{lem:LaplacianHarmonicFunctionPower}
\cite[Equation (3.10)]{GroisserParkerDecay}
Let $(X,g)$ be a Riemannian, smooth $m$-manifold, let $h$ be a smooth 
positive harmonic function on $X$, let $\xi = |d(\log\sqrt{h})|^2$,
and let $\alpha$ be a constant. Then
\begin{equation}
\label{eq:LaplacianHarmonicFunctionPower}
\Delta h^\alpha
=
4\alpha(1-\alpha)\xi h^\alpha.
\end{equation}
\end{lem}

\begin{proof}
One simply calculates that
$$
\Delta h^\alpha
=
d^*dh^\alpha
=
\alpha(1-\alpha)h^\alpha|h^{-1}dh|^2
=
4\alpha(1-\alpha)\xi h^\alpha,
$$
as desired.
\end{proof}

We apply Lemma \ref{lem:LaplacianHarmonicFunctionPower} to the harmonic
function $h$ on $(X,g)$. We observe that
$$
4\alpha(1-\alpha) = \delta(2-\delta),
$$
for $\alpha = \delta/2$ or $\alpha = 1 - \delta/2$. Hence, Lemma
\ref{lem:LaplacianHarmonicFunctionPower} yields
\begin{equation}
\label{eq:RiemLaplacianHarmFun2TestPowers}
\begin{aligned}
\Delta h^{\delta/2} 
&= \delta(2-\delta)\xi h^{\delta/2},
\\
\Delta h^{1-\delta/2}
&= \delta(2-\delta)\xi h^{1-\delta/2}.
\end{aligned}
\end{equation}
Therefore, we define
$$
g_1
=
r_0^{-\delta+2/3}h^{1-\delta/2} 
+ r_1^{\delta-4/3}h^{\delta/2}.
$$
Then the definition of $g_1$, the fact that $0<\delta <1$, the
inequalities \eqref{eq:HarmonicFunction}, and the identities in
\eqref{eq:RiemLaplacianHarmFun2TestPowers} imply that
\begin{equation}
\label{eq:DiffIneqLaplaceg1}
\begin{cases}
g_1 \geq \half r_0^{-4/3} &\text{on }\partial B(x_0,2r_0),
\\
g_1 \geq \half r_1^{-4/3} &\text{on }\partial B(x_0,r_1/2),
\\
(\Delta-\delta(2-\delta)\xi)g_1
= 0 &\text{on } \Omega(x_0,2r_0,r_1/2).
\end{cases}
\end{equation}
Thus, it follows from Lemma \ref{lem:BasicComparison} that, on
$\Omega(x_0,2r_0,r_1/2)$,
\begin{equation}
\label{eq:PrelimPhiDecayr}
\begin{aligned}
|\phi|^{2/3} \leq cE^{2/3}g_1 
&=
cE^{2/3}\left(r_0^{-\delta+2/3}h^{1-\delta/2} 
+ r_1^{\delta-4/3}h^{\delta/2}\right)
\\
&=
cE^{2/3}\left(r_0^{-4/3}(r_0^2h)^{1-\delta/2} 
+ r_1^{-4/3}(r_1^2h)^{\delta/2}\right).
\end{aligned}
\end{equation}
If we now substitute our hypothesized decay estimates
\eqref{eq:ScalarCurvature} for $|\kappa_g|$ and the estimate for $|F_A|$ on
$\Omega(x_0,r_0,r_1)$ given by \eqref{eq:RadeDecay},
\eqref{eq:GPAsympFlatDecay}, or \eqref{eq:GPPuncturedBallDecay}) (courtesy
of \eqref{eq:SpincCurvature}), and the preliminary decay estimate
\eqref{eq:PrelimPhiDecayr} for $\phi$ into the differential inequality
\eqref{eq:RawHarmonicSectionTwoThirdDiffIneq} for $\Delta|\phi|^{2/3}$, we
obtain (using $r_0\geq 1$, $2r_0<r<r_1/2$, $\delta < 1$, and inequalities
\eqref{eq:HarmonicFunction} for $h$)
\begin{equation} 
\label{eq:RiemQuarticLaplacePhiTwoThird}
\begin{aligned}
\Delta|\phi|^{2/3}
&\leq
c\eps(r_0^2r^{-4}+r_1^{-2})|\phi|^{2/3} \quad\text{on }\Omega(x_0,2r_0,r_1/2)
\\
&\leq
c\eps\left(r_0^{-4/3}(r_0^2h)^{3-\delta/2} 
+ r_1^{-4/3}(r_1^2h)^{\delta/2}\right).
\end{aligned}
\end{equation} 
Therefore, substituting the preliminary decay estimate
\eqref{eq:PrelimPhiDecayr} for $\phi$ into
\eqref{eq:RiemQuarticLaplacePhiTwoThird}, together with the upper bound
\eqref{eq:HarmonicFunction} for $r^{-2}$, yields
\begin{equation}
\label{eq:BetterDiffIneqLaplacePhiTwoThird}
\begin{cases}
|\phi|^{2/3} \leq cE^{2/3}r_0^{-4/3} &\text{on }\partial B(x_0,2r_0),
\\
|\phi|^{2/3} \leq cE^{2/3}r_1^{-4/3} &\text{on }\partial B(x_0,r_1/2),
\\
\Delta|\phi|^{2/3}
\leq c\eps E^{2/3}
\left(r_0^{-4/3}(r_0^2h)^{3-\delta/2} + r_1^{-4/3}(r_1^2h)^{\delta/2}\right)
&\text{on } \Omega(x_0,2r_0,r_1/2).
\end{cases}
\end{equation}
Define
\begin{align*}
g_2
=
2\left(r_0^{2/3}h + r_1^{-4/3}\right)
-
\left(r_0^{-4/3}(r_0^2h)^{2-\delta/2} + r_1^{-4/3}(r_1^2h)^{\delta/2-1}\right),
\end{align*}
and note that
\begin{equation}
\label{eq:g2UpperBound}
g_2
\leq
2\left(r_0^{2/3}h + r_1^{-4/3}\right),
\end{equation}
while, using the inequalities \eqref{eq:HarmonicFunction} for $h$, the fact
that $2r_0<r<r_1/2$, and $\delta<1$, gives
\begin{equation}
\label{eq:g2LowerBound}
g_2
\geq
\left(r_0^{2/3}h + r_1^{-4/3}\right)
\quad\text{on }\Omega(x_0,2r_0,r_1/2).
\end{equation}
Observe that the identity \eqref{eq:LaplacianHarmonicFunctionPower} yields
\begin{align*}
\Delta h^{2-\delta/2} &= -(4-\delta)(2-\delta)\xi h^{2-\delta/2},
\\
\Delta h^{\delta/2-1} &= -(4-\delta)(2-\delta)\xi h^{\delta/2-1}.
\end{align*}
Then, by definition of $g_2$, the lower bound
\eqref{eq:HarmonicFunction} for $\xi$, and the preceding inequalities
on $\partial\Omega(x_0,2r_0,r_1/2)$, we see that
\begin{equation}
\label{eq:DiffIneqLaplaceg2}
\begin{cases}
g_2 \geq \half r_0^{-4/3} &\text{on }\partial B(x_0,2r_0),
\\
g_2 \geq \half r_1^{-4/3} &\text{on }\partial B(x_0,r_1/2),
\\
\begin{aligned}
\Delta g_2 &\geq c(4-\delta)(2-\delta)
\\
&\quad\times\left(r_0^{-4/3}(r_0^2h)^{3-\delta/2} +
r_1^{-4/3}(r_1^2h)^{\delta/2}\right)
\end{aligned}
&\text{on } \Omega(x_0,2r_0,r_1/2).
\end{cases}
\end{equation}
Note that the constant $(4-\delta)(2-\delta)$ is positive as we assumed
$\delta < 1$. Therefore, it follows from the inequalities
\eqref{eq:BetterDiffIneqLaplacePhiTwoThird} and
\eqref{eq:DiffIneqLaplaceg2} that
\begin{align*}
|\phi|^{2/3} &\leq cE^{2/3}g_2\text{ on }\partial\Omega(x_0,2r_0,r_1/2),
\\
\Delta |\phi|^{2/3} &\leq \Delta(c\eps E^{2/3}g_2) \leq
\Delta(cE^{2/3}g_2) \text{ on }\Omega(x_0,2r_0,r_1/2). 
\end{align*}
Hence, the comparison principle for the Laplacian $\Delta$ and the upper
bound \eqref{eq:g2UpperBound} for $g_2$ that on $\Omega(x_0,2r_0,r_1/2)$
implies that
$$
|\phi|^{2/3} \leq cE^{2/3}g_2 
\leq
cE^{2/3}\left(r_0^{2/3}h + r_1^{-4/3}\right)
$$
and so
$$
|\phi|
\leq cE\left(r_0h^{3/2} + r_1^{-2}\right)
\leq cE\left(r_0r^{-3} + r_1^{-2}\right).
$$
This completes the proof of Theorem \ref{thm:Main} when $k=0$. The cases
$k\geq 1$ follow from immediately from the case $k=0$ and Lemma
\ref{lem:BoundaryDecayEstimate}.

\subsection{Decay estimates for eigenspinors
and proof of Corollary \ref{cor:SEMain}}
\label{subsec:EigenspinorDecay}
We now allow $\phi$ to be an eigenspinor of $D_A$, with non-zero eigenvalue
$\mu$. Thus $\phi$ is a $D_{\tilde A}$-harmonic spinor when $\tilde A$ is
the Friedrich connection \eqref{eq:FriedrichConnDim4} defined by $A$ and
$\mu$. The curvature $F_{\tilde A}$ obeys
$$
|F_{\tilde A}| \leq |F_A| + c\mu^2,
$$
where $c$ is a universal constant (independent of the metric). We must
assume that $\mu^2$ satisfies the inequalities \eqref{eq:ScalarCurvature}
obeyed by $\kappa_g$ for suitable $r_0\leq r\leq r_1$, so in the argument
of \S \ref{subsec:HarmonicSpinorDecay} the bounds
\eqref{eq:RiemQuadraticLaplacePhiTwoThird} and
\eqref{eq:RiemQuarticLaplacePhiTwoThird} continue to hold with $|\kappa_g|$
replaced by $|\kappa_g|+\mu^2$ when $\mu\neq 0$. Writing
$\mu^2=\mu^2r_1^2r_1^{-2}$, one trivially has
$$
\mu^2 \leq \eps r_1^{-2} \leq \eps r^{-2}, \quad r_0\leq r\leq r_1,
$$
provided $r_1$ is small enough that $r_1^2\leq\eps$.  Then the inequalities
\eqref{eq:ScalarCurvature} hold with $\mu^2$ in place of $\kappa_g$ when
$c_1=1$ and $c_2=r_1^2$. Hence, the same argument which proves Theorem
\ref{thm:Main} leads to a proof of Corollary \ref{cor:SEMain}, since the
additional term $F_{\tilde A}-F_A$ obeys the same estimates (when
$r_1^2\leq\eps$) and has the same scaling behavior as $\kappa_g$.

\section{Bubbling and elliptic estimates for the Dirac operator}
\label{sec:EllipticDirac}
Theorem \ref{thm:Main} and Corollary \ref{cor:SEMain} provide pointwise
decay estimates for eigenspinors which are useful over annuli in a
four-manifold $(X,g)$ where we have suitable bounds on the curvatures of
$A$ and the Levi-Civita connection $\cov_g$. In this section we describe
(see Proposition \ref{prop:MainNegEstimate}) how one can apply Theorem
\ref{thm:Main} to the problem of obtaining $C^0\cap L^2_2$ estimates for
negative harmonic spinors when the induced $\SO(3)$ 
connection $A$ on $\su(E)$ ``bubbles'' in the
Uhlenbeck sense \cite{UhlLp}. Estimates of this kind have important
applications to the problem of gluing $\PU(2)$ monopoles \cite{FL3},
\cite{FL4}. We emphasize negative spinors because it is the
anti-self-dual component of the curvature which appears in the Bochner
formula \eqref{eq:Bochner} for $D_A^2$ on negative spinors, whereas the
self-dual component of the curvature appears in the Bochner for $D_A^2$ on
positive spinors: with the usual conventions of gauge theory (for the
anti-self-dual Yang-Mills equation \cite{DK} or the $\PU(2)$ monopole equations
\cite{FL1}), the anti-self-dual component of the curvature can
bubble while one has good control on the self-dual component. Some remarks
on the significance of Proposition \ref{prop:MainNegEstimate} are included
in \S \ref{subsec:RemarksOnWhyPropIsInteresting}.

\subsection{Integral estimates for negative spinors when the anti-self-dual
curvature component bubbles} To describe our application, we continue with
the setup of \S \ref{subsec:Statement} and write $V=V^+\oplus V^-$, where
$V^\pm = W^\pm\otimes E$. We further suppose that $E$ is complex rank two,
let $\hat A$ denote the induced $\SO(3)$ connection on $\su(E)$, and let
$A_e$ denote the induced unitary connection on the line bundle
$\det(E)$. In gauge-theoretic applications \cite{DK}, \cite{FL1}, it is the
curvature component $F_{\hat A}$ which may bubble, while the remaining
connections (on $W$ and $\det(E)$) are usually fixed.  We allow the
constant $c$ in our estimates below to depend on $C^2$ bounds on the
curvature of the Riemannian metric $g$ and the connections $A_d$ and $A_e$.

Recall \cite[\S 4]{FeehanSliceV2}, \cite{TauStable} that the $L^{\sharp}$
and $L^{2\sharp}$ Sobolev norms are defined by
\begin{align*}
\|\phi\|_{L^{\sharp}(X)}
:=
\sup_{x\in X}\|\dist^{-2}(x,\cdot)|\phi|\|_{L^1(X)}
&\quad\text{and}\quad
\|\phi\|_{L^{\sharp,2}(X)}
:= 
\|\phi\|_{L^{\sharp}(X)} + \|\phi\|_{L^2(X)}.
\\
\|\phi\|_{L^{2\sharp}(X)}
:=
\sup_{x\in X}\|\dist^{-1}(x,\cdot)|\phi|\|_{L^2(X)}
&\quad\text{and}\quad
\|\phi\|_{L^{2\sharp,4}(X)}
:= 
\|\phi\|_{L^{2\sharp}(X)} + \|\phi\|_{L^4(X)}.
\end{align*}
See \cite[\S 4]{FeehanSliceV2}, \cite{TauStable} for further explanation and
properties of this family of Sobolev norms.  

\begin{prop}
\label{prop:MainNegEstimate}
Let $X$ be a closed, oriented four-manifold with metric $g$. Given $m\in
\NN$, there are positive constants $r_1$, $c(m,r_1)$, and $\eps$ such that
for small enough $\eps$ the following holds.  Suppose $\{x_i\}_{i=1}^m$
is a finite set of distinct points and $\{\lambda_i\}_{i=1}^m\subset
(0,r_1]$ is a set of positive constants. Let
$$
U  = X-\bigcup_{i=1}^m B(x_i,{\textstyle{\frac{1}{2}}}\lambda_i^{1/2})
\quad\text{and}\quad
U' = X-\bigcup_{i=1}^m B(x_i,4\lambda_i^{1/3})\Subset U.
$$
Assume the $L^2$ norms of the curvatures of $A_d$, $A_e$, and $\hat A$
on the annuli $\Omega(x_i,\frac{1}{2}\lambda_i^{1/3},r_1)$ 
are less than or equal to $\eps$. If $\phi \in C^\8(X,V^-)\cap \Ker D_A$, then
\begin{align*}
\|\phi\|_{C^0\cap L^2_{2,A}(U')}
&\le
c(1+\|F_{\hat A}\|_{L^2(X)})(1+\|F_{\hat A}^-\|_{C^0(U)})^2\|\phi\|_{L^2(U)}.
\end{align*}
\end{prop}

In a typical application of of Proposition \ref{prop:MainNegEstimate}, the
curvature of $\hat A$ would be allowed to bubble near the points $x_i\in
X$. Before proceeding with the proof, we recall the following important
estimate from \cite{FeehanSliceV2}:

\begin{lem}
\label{lem:L22InfinityEstu}
\cite[Lemma 5.5]{FeehanSliceV2}
Let $X$ be a closed, oriented, Riemannian four-manifold. Then there is
a constant $c$ with the following significance. Let $\sF$ be a Riemannian
vector bundle over $X$ and let $B$ be an orthogonal
$C^\8$ connection on $\sF$ with curvature $F_{B}$. Then the following
estimate holds for any $\phi\in C^\8(\sF)$: 
\begin{align*}
\|\phi\|_{C^0\cap L^2_{2,B}(X)}
&\le c(1+\|F_{B}\|_{L^2(X)})
(\|\cov_{B}^*\cov_{B}\phi\|_{L^{\sharp,2}(X)} + \|\phi\|_{L^2(X)}). 
\end{align*}
\end{lem}

\begin{proof}[Proof of Proposition \ref{prop:MainNegEstimate}]
We choose a smooth cutoff function $\beta$ so that (see the proof of
Lemma 5.8 in \cite{FLKM1V3})
\begin{equation}
\label{eq:CutoffFnDerivBound}
\|\cov\beta\|_{L^{2\sharp,4}(X)} +  \|\cov^2\beta\|_{L^{\sharp,2}(X)}\leq c,
\end{equation}
for some constant $c(g)$, where
$$
\beta
= 
\begin{cases} 
1 &\text{on }X-\cup_{i=1}^m B(x_i,2\lambda_i^{1/3}), 
\\
0 &\text{on }\cup_{i=1}^m B(x_i,\lambda_i^{1/3}).
\end{cases}
$$
For the bound on $\|\cov^2\beta\|_{L^{\sharp}(X)}$ one uses the fact that
this norm is scale invariant (it has the same scaling behavior as
$\|\cov^2\beta\|_{L^{2}(X)}$), so we can assume that $\lambda_i=1$ without
loss of generality when deriving the uniform bound on
$\|\cov^2\beta\|_{L^{\sharp}(X)}$. 

Note that $\beta=1$ on $U'$ and $\supp\beta\subset U$. Schematically, we have
$$
\cov_A^*\cov_A(\beta\phi)
=
(\Delta\beta)\phi + \cov\beta\otimes\cov_A\phi
+ \beta\cov_A^*\cov_A\phi,
$$
and therefore the cutoff function derivative bounds
\eqref{eq:CutoffFnDerivBound} yield
\begin{equation}
\label{eq:CovLapacianCutoffSpinor}
\begin{aligned}
\|\cov_A^*\cov_A(\beta\phi)\|_{L^{\sharp,2}(X)}
&\leq
c\|\phi\|_{C^0(\supp d\beta)}
+ c\|\cov_A\phi\|_{L^{2\sharp,4}(\supp d\beta)}
\\
&\quad+ \|\cov_A^*\cov_A\phi\|_{L^{\sharp,2}(\supp\beta)}.
\end{aligned}
\end{equation}
Since $D_A^2\phi=0$, the Bochner formula \eqref{eq:Bochner} implies that
\begin{equation}
\label{eq:BochnerCovLapSpinorEstimate}
\|\cov_A^*\cov_A\phi\|_{L^{\sharp,2}(\supp\beta)}
\leq 
c(1+\|F_{\hat A}^-\|_{C^0(\supp\beta)})\|\phi\|_{L^{\sharp,2}(\supp\beta)}.
\end{equation}
Theorem \ref{thm:Main} gives a $C^0$ bound for
$\phi$ on $\supp d\beta \subset
\cup_{i=1}^m\Omega(x_i,\lambda_i^{1/3},2\lambda_i^{1/3})$:
\begin{equation}
\label{eq:C0EstSpinorOnAnnuli}
\|\phi\|_{C^0(\supp d\beta)}
\leq
\sum_{i=1}^m\|\phi\|_{C^0(\Omega(x_i,\lambda_i^{1/3},2\lambda_i^{1/3}))}
\leq
cm\|\phi\|_{L^2(U)},
\end{equation}
where we take $r_0=\lambda_i$. Lemma 4.1 in \cite{FeehanSliceV2} implies that
$\|\phi\|_{L^{\sharp}(\supp\beta)}\leq c\|\phi\|_{L^4(\supp\beta)}$. We
may consider $U$ as the union of the region
$U_0 = X-\cup_{i=1}^mB(x_i,\frac{1}{2}r_1)$ and the annuli
$\Omega(x_i,\lambda_i^{1/3},r_1)$. We can choose a
cutoff function $\chi$ which is equal to one $U_0$ and supported on the
complement of the balls $B(x_i,\frac{1}{4}r_1)$. We now apply the
Sobolev embedding $L^2_{1,A}\subset L^4$ (with embedding constant
independent of $A$), integration by 
parts, and the Bochner formula \eqref{eq:Bochner} to obtain
$$
\|\chi\phi\|_{L^4(X)}
\leq
c\|\chi\phi\|_{L^2_{1,A}(X)}
\leq
c(1+\|F_{\hat A}^-\|_{C^0(U_0)})\|\chi\phi\|_{L^2(X)}.
$$
Over the annuli $\Omega(x_i,\lambda_i^{1/3},r_1)$, Theorem \ref{thm:Main}
yields bounds of the form
$$
\|\phi\|_{L^2_{1,A}(\Omega(x_i,\lambda_i^{1/3},r_1))}
\leq 
c\|\phi\|_{L^2(U)}.
$$
Hence, combining these inequalities, we have
\begin{equation}
\label{eq:phiEstLSharpOnSuppBeta}
\|\phi\|_{L^{\sharp}(\supp\beta)}
\leq 
c(1+\|F_{\hat A}^-\|_{C^0(U_0)})\|\phi\|_{L^2(U)}.
\end{equation}
We next consider the term $\|\cov_A\phi\|_{L^{4}(\supp d\beta)}$ on the
right-hand side of estimate \eqref{eq:CovLapacianCutoffSpinor}.
Denoting $C_1=\|\phi\|_{L^2(U)}$, we observe that Theorem \ref{thm:Main}
(again taking $r_0=\lambda_i$) gives
\begin{align*}
\|\cov_A\phi\|_{L^4(\supp d\beta)}^4
&=
cC_1^4\sum_{i=1}^m\int_{\lambda_i^{1/3}}^{2\lambda_i^{1/3}}
r^{-4}(\lambda_i r^{-3} + r_1^2)^4 r^3\,dr 
\\
&\leq cC_1^4\sum_{i=1}^m\int_{\lambda_i^{1/3}}^{2\lambda_i^{1/3}} r^{-1}\,dr
\\
&= 
\left. cC_1^4\sum_{i=1}^m \log r\right|_{\lambda_i^{1/3}}^{2\lambda_i^{1/3}}
= cmC_1^4\log 2,
\end{align*}
and therefore
\begin{equation}
\label{eq:L4CovSpinorBound}
\|\cov_A\phi\|_{L^4(\supp d\beta)}
\leq
c\|\phi\|_{L^2(U)}.
\end{equation}
For the term $\|\cov_A\phi\|_{L^{2\sharp}(\supp d\beta)}$ on the right-hand
side of estimate \eqref{eq:CovLapacianCutoffSpinor}, choose a cutoff
function $\gamma$ defined in much 
the same way as $\beta$, except that $\gamma=0$ on the complement of the
annulus $\Omega(x_i,\frac{1}{2}\lambda_i^{1/3},4\lambda_i^{1/3})$ and
$\gamma=1$ on the annulus
$\Omega(x_i,\lambda_i^{1/3},2\lambda_i^{1/3})$. Thus, we 
have $\supp\gamma\subset U$. Then
\begin{equation}
\label{eq:L2SharpBoundOnCovPhiAnnuli}
\begin{aligned}
\|\cov_A\phi\|_{L^{2\sharp}(\supp d\beta)}
&\leq
\|\gamma\cov_A\phi\|_{L^{2\sharp}(X)}
\leq
\|\gamma\cov_A\phi\|_{L^{2}_{1,A}(X)}
\quad\text{(Lemma 4.1 in \cite{FeehanSliceV2})}
\\
&\leq
c\|\gamma\cov_A\phi\|_{L^{4}(X)}
+ \|\cov_A(\gamma\cov_A\phi)\|_{L^{2}(X)}
\\
&\leq
c\left(\|\cov_A\phi\|_{L^{4}(\supp\gamma)}
+ \|\cov_A^2\phi\|_{L^{2}(\supp\gamma)}\right).
\end{aligned}
\end{equation}
The final inequality above follows from H\"older's inequality and the
uniform $L^4$ bound on $\cov\gamma$, analogous to that in
\eqref{eq:CutoffFnDerivBound} for $\cov\beta$. The argument yielding inequality
\eqref{eq:L4CovSpinorBound} also provides 
a similar bound on $\|\cov_A\phi\|_{L^{4}(\supp\gamma)}$.
For the term $\|\cov_A^2\phi\|_{L^{2}(\supp\gamma)}$, Theorem
\ref{thm:Main} gives
\begin{align*}
\|\cov_A^2\phi\|_{L^2(\supp\gamma)}^2
&=
cC_1^2\sum_{i=1}^m\int_{\frac{1}{2}\lambda_i^{1/3}}^{4\lambda_i^{1/3}}
r^{-4}(\lambda_i r^{-3} + r_1^2)^2 r^3\,dr 
\\
&\leq
cC_1^2\sum_{i=1}^m\int_{\frac{1}{2}\lambda_i^{1/3}}^{4\lambda_i^{1/3}} 
r^{-1}\,dr
\\
&= 
cmC_1^2\log 2,
\end{align*}
and therefore, combining estimates \eqref{eq:L4CovSpinorBound} 
(now for $\|\cov_A\phi\|_{L^{4}(\supp\gamma)}$) and
\eqref{eq:L2SharpBoundOnCovPhiAnnuli} with the preceding inequality, we obtain
\begin{equation}
\label{eq:L2SharpCovSpinorBound}
\|\cov_A\phi\|_{L^{2\sharp}(\supp d\beta)}
\leq
c\|\phi\|_{L^2(U)}.
\end{equation}
Combining the inequalities \eqref{eq:CovLapacianCutoffSpinor},
\eqref{eq:BochnerCovLapSpinorEstimate}, 
\eqref{eq:C0EstSpinorOnAnnuli}, \eqref{eq:phiEstLSharpOnSuppBeta},
\eqref{eq:L4CovSpinorBound}, and 
\eqref{eq:L2SharpCovSpinorBound} gives 
$$
\|\cov_A^*\cov_A(\beta\phi)\|_{L^{\sharp,2}(X)}
\leq
c(1+\|F_{\hat A}^-\|_{C^0(U)})^2\|\phi\|_{L^{2}(U)}.
$$
Substituting the preceding bound into the estimate given by Lemma
\ref{lem:L22InfinityEstu} for $\beta\phi$ yields
\begin{align*}
\|\phi\|_{C^0\cap L^2_{2,A}(U')}
&\leq
\|\beta\phi\|_{C^0\cap L^2_{2,A}(X)}
\\
&\leq
c(1+\|F_{\hat A}\|_{L^2(X)})
(1+\|F_{\hat A}^-\|_{C^0(U)})^2\|\phi\|_{L^{2}(U)}.
\end{align*}
This completes the proof.
\end{proof}

\subsection{Remarks on the result}
\label{subsec:RemarksOnWhyPropIsInteresting}
In contrast to the proof of Proposition \ref{prop:MainNegEstimate}, which
relied on our pointwise decay estimates for negative harmonic spinors, the
Bochner formula \eqref{eq:Bochner} and a simple integration-by-parts
argument yields a useful estimate for positive harmonic spinors:

\begin{lem}
\label{lem:LinftyL22Estphi}
Let $X$ be a closed, oriented, Riemannian four-manifold. Then there are
positive constants $c$, $\eps$ so that for small enough $\eps$ the
following holds. Suppose $X=U_2\cup U_\8$, where $U_2\Subset X$
is an open subset with $\|F_{\hat A}^+\|_{L^{\sharp,2}(U_2)}
<\eps$. If $\phi\in C^\8(X,V^+)\cap \Ker D_A$, then
\begin{align*}
\|\phi\|_{C^0\cap L^2_{2,A}(X)}
&\leq
c(1+\|F_{\hat A}\|_{L^2(X)})^2(1+\|F_{\hat A}^+\|_{C^0(U_\8)})^2
\|\phi\|_{L^2(X)}. 
\end{align*}
\end{lem}

\begin{proof}
{}From the Bochner formula \eqref{eq:Bochner} and the Sobolev multiplication
\cite[Lemma 4.3]{FeehanSliceV2},
\begin{align*}
\|\cov_A^*\cov_A\phi\|_{L^{\sharp,2}}
\le c(1+\|F_{\hat A}^+\|_{C^0(U_\8)})\|\phi\|_{L^{\sharp,2}}
+ c\|F_{\hat A}^+\|_{L^{\sharp,2}(U_2)}\|\phi\|_{C^0(U_2)}.
\end{align*}
Combining the preceding estimate with that of Lemma
\ref{lem:L22InfinityEstu}, together with the embedding and interpolation
inequalities $\|\phi\|_{L^\sharp} \le c\|\phi\|_{L^4}\le
c\|\phi\|_{L^2}^{1/2}\|\phi\|_{C^0}^{1/2}$, choosing $\eps$ small enough so
that $c\|F_{\hat A}^+\|_{L^{\sharp,2}(U_2)}\leq \frac{1}{2}$, and using
rearrangement with the terms involving $\|\phi\|_{C^0}$ then yields the
desired bound.
\end{proof}

Lemma \ref{lem:LinftyL22Estphi} provides useful $C^0\cap L^2_2(X)$ elliptic
estimates for positive spinors $\phi\in C^\8(X,V^+)$ even when $\hat A$
bubbles, because we still have uniform $C^0$ bounds for $F^+_{\hat A}$ away
from the bubble points (on the set $U_\8$) and small $L^{\sharp,2}$ bounds
around the bubble points (on the set $U_2$). However, this is never the
case for $F^-_{\hat A}$ in such applications because $F^-_{\hat A}$ is
neither $C^0$-bounded nor $L^{\sharp,2}$-small around the bubble points. So
the analogue of Lemma \ref{lem:LinftyL22Estphi} for negative
spinors---which would have $F^-_{\hat A}$ in place of $F^+_{\hat 
A}$---does not provide useful $C^0\cap L^2_2(X)$ elliptic estimates for
$\phi\in C^\8(X,V^-)$. In contrast, in Proposition
\ref{prop:MainNegEstimate}---where $U$ plays the role of $U_\8$ above---we
at least still obtain an estimate for $\phi$ on the complement of the small
balls $B(x_i,4\lambda_i^{1/3})$, which is uniform with respect $\lambda_i$
as those constants tend to zero.


\ifx\undefined\bysame
\newcommand{\bysame}{\leavevmode\hbox to3em{\hrulefill}\,}
\fi


\end{document}